\documentclass{article}%
\usepackage{amssymb}
\usepackage{amsmath}
\usepackage{amsfonts}
\usepackage{graphicx}%
\providecommand{\U}[1]{\protect\rule{.1in}{.1in}}

\begin{document}

\title{What does the proof of Birnbaum's theorem prove?}
\author{\ Michael Evans\\Department of Statistics\\University of Toronto}
\date{}
\maketitle

\noindent\textit{Abstract}: Birnbaum's theorem, that the sufficiency and
conditionality principles entail the likelihood principle, has engendered a
great deal of controversy and discussion since the publication of the result
in 1962. In particular, many have raised doubts as to the validity of this
result. Typically these doubts are concerned with the validity of the
principles of sufficiency and conditionality as expressed by Birnbaum.
Technically it would seem, however, that the proof itself is sound. In this
paper we use set theory to formalize the context in which the result is proved
and show that in fact Birnbaum's theorem is incorrectly stated as a key
hypothesis is left out of the statement. When this hypothesis is added, we see
that sufficiency is irrelevant, and that the result is dependent on a
well-known flaw in conditionality that renders the result almost
vacuous.\bigskip


\section{Introduction}

A result presented in Birnbaum (1962), and referred to as Birnbaum's theorem,
is very well-known in statistics. This result says that a statistician who
accepts both the sufficiency $S$ and conditionality $C$ principles must also
accept the likelihood principle $L$ and conversely. The result has always been
controversial primarily because it implies that a frequentist statistician who
accepts $S$ and $C$ is forced to ignore the repeated sampling properties of
any inferential procedures they use. Given that both $S$ and $C$ seem quite
natural to many frequentist statisticians while $L$ does not, the result is
highly paradoxical.

Various concerns have been raised about the proof of the result. For example,
Durbin (1970) argued that the theorem fails to hold whenever $C$ is restricted
by requiring that any ancillaries used must be functions of a minimal
sufficient statistic. Kalbfleisch (1975) argued that $C$ should only be
applicable when the value of the ancillary statistic used to condition is
actually a part of the experimental make-up. This is called the weak
conditionality principle. In Evans, Fraser and Monette (1986) it is argued
that Birnbaum's theorem, and a similar result that accepting $C$\ alone is
equivalent to accepting $L$, are invalid because the specific uses of $S$ and
$C$ in proving these results can be seen to be based on flaws in their
formulations. For example, Birnbaum's theorem requires a use of $S$ and $C$
where the information discarded by $S$ as irrelevant, which is the primary
motivation for $S,$ is exactly the information used by $C$ to condition on and
so identifies the discarded information as highly relevant. As such $S$ and
$C$ contradict each other. We note that this is precisely what Durbin's
restriction on the ancillaries avoids. Furthermore, the result that $C$ alone
implies $L$ can be seen to depend on the lack of a unique maximal ancillary
which can be viewed as an essential flaw in $C$. Also, see Holm (1985),
Barndorff-Nielsen (1995) and Helland (1995) for various concerns about the
formulation of the theorem. Mayo (2010) argues that, in the context of a
repeated sampling formulation for statistics, we cannot simultaneously have
$S$ and $C$ true, as when $S$ is true then $C$ is false and when $C$ is true
then $S$ is false. Gandenberger (2012) offers up a proof that avoids some of
the objections raised by others.

Many of these reservations are essentially with the hypotheses to the theorem
and suggest that Birnbaum's theorem should be rejected because the hypotheses
are either not acceptable or have been misapplied. It is the purpose of this
paper to provide a careful set-theoretic formulation of the context of the
theorem. When this is done we see that there is a hypothesis that needs to be
formally acknowledged as part of the statement of Birnbaum's theorem. With
this addition, the force of the result is lost and the paradox disappears. The
same conclusions apply to result that $C$ is equivalent to $L$ and, in fact,
this is really the only result as $S$ is redundant in Birnbaum's theorem when
the additional hypothesis is formally acknowledged.

For our discussion it is important that we stick as closely as possible to
Birnbaum's formulation. To discuss the proof, however, we have to make certain
aspects of Birnbaum's argument mathematically precise that are somewhat vague
in his paper. It is always possible then that someone will argue that we have
done this in a way that is not true to Birnbaum's intention.\ We note,
however, that this is accomplished in a very simple and direct way. If there
is another precise formulation that makes the theorem true, then it is
necessary for a critic of how we do this to provide that alternative.

A basic step missing in Birnbaum (1962) was to formulate the principles as
relations on the set $\mathcal{I}$ of all model and data combinations. So
$\mathcal{I}$ is the set of all \textit{inference bases} $I=(E,x)$ where
$E=(\mathcal{X}_{E},\{f_{E,\theta}:\theta\in\Theta_{E}\}),$ $\mathcal{X}_{E}$
is a sample space, $\{f_{E,\theta}:\theta\in\Theta_{E}\}$ is a collection of
probability density functions on $\mathcal{X}_{E},$ with respect to some
support measure $\mu_{E}$ on $\mathcal{X}_{E},$ and $x\in\mathcal{X}_{E}$ is
the observed data. We will ignore all measure-theoretic considerations as they
are not essential for any of the arguments. If the reader is concerned by
this, then we note that the collection of models where $\mathcal{X}_{E}$ and
$\Theta_{E}$ are finite and $\mu_{E}$ is counting measure is rich enough to
produce the paradoxical result. So in general we can consider our discussion
restricted to the case where $\mathcal{X}_{E}$ and $\Theta_{E}$ are finite. It
is our view that infinite sets and continuous probability measures are not
necessary for the development of the basic principles of statistics. Rather
the use of infinite sets and continuity represents approximations to a finite
reality and appropriate restrictions must be employed on such quantities so
that we are not mislead by purely mathematical considerations. In spite of our
restrictions, most of our development applies equally well under very general circumstances.

We note that expressing the principles as relations was part of Evans, Fraser
and Monette (1986) this is taken further here. In Section 2 we discuss the
meaning and use of relations generally. In Section 3 we apply our discussion
of relations to Birnbaum's theorem. In Section 4 we draw some conclusions.

\section{Relations}

A \textit{relation }$R$ with domain $D$ is a subset $R\subset D\times D.$
Saying $(x,y)\in R$ means that the objects $x$ and $y$ have a property in
common. For example, suppose $D$ is the set of students enrolled at a specific
university at a specific point in time. Let $R_{1}$ be defined by $(x,y)\in
R_{1}$ whenever $x$ and $y$ are students in the same class. Let $R_{2}$ be
defined by $(x,y)\in R_{2}$ whenever $x$ and $y$ have taken a course from the
same professor.

A relation $R$ is \textit{reflexive} if $(x,x)\in R$ for all $x\in D$,
\textit{symmetric} if $(x,y)\in R$ implies $(y,x)\in R,$ and
\textit{transitive} if $(x,y)\in R,(y,z)\in R$ implies that $(x,z)\in R.$ If a
relation $R$ is reflexive, symmetric and transitive, then $R$ is called an
\textit{equivalence relation}. Clearly $R_{1}$ is an equivalence relation and,
while $R_{2}$ is reflexive and symmetric, it is not typically transitive and
so is not an equivalence relation. While $(x,y)\in R$ implies that $x$ and $y$
are related, perhaps by the possession of some property, when $R$ is an
equivalence relation this implies that $x$ and $y$ possess the property to the
same degree. We say that relation $R$ on $D$ \textit{implies} relation
$R^{\prime}$ on $D$ whenever $R\subset R^{\prime}.$ Clearly we have \ that
$R_{1}\subset R_{2}.$

If $R$ is a relation on $D$, then the equivalence relation $\bar{R}$ generated
by $R$ is the smallest equivalence relation containing $R.$ We see that
$\bar{R}$ is the intersection of all equivalence relations on $D$ containing
$R.$ Also we have that
\begin{align}
\bar{R}=\{  &  (x,y):\exists\,n,x_{1},\ldots,x_{n}\in D\text{ with }%
x=x_{1},y=x_{n}\text{ and }\nonumber\\
&  (x_{i},x_{i+1})\in R\text{ or }(x_{i+1},x_{i})\in R\}. \tag{1}%
\end{align}

It is not always clear that $\bar{R}$ has a meaningful interpretation, at
least as it relates to the property being expressed by $R.$ For example,
$\bar{R}_{2}$ is somewhat more difficult to interpret and surely goes beyond
the idea that $R_{2}$ is perhaps trying to express, namely, that two students
were directly influenced by the same professor. In fact, it is entirely
possible that $\bar{R}_{2}=D\times D.$ As another example, suppose that
$D=\{2,3,4,\ldots\}$ and $(x,y)\in R$ when $x$ and $y$ have a common factor
bigger than 1. Then $R$ is reflexive and symmetric but not transitive. If
$x,y\in D$ then $(x,xy)\in R,(xy,y)\in R$ so $\bar{R}=D\times D$ and $\bar{R}$
is saying nothing. It seems that each situation, where we extend a relation
$R$ to an equivalence relation, must be examined to see whether or not this
extension has any meaningful content for the application.

Now suppose we have relations $R_{1}$ and $R_{2}$ on $D$ and consider the
relation $R_{1}\cup R_{2}.$ The following result is relevant to our discussion
in Section 3.$\smallskip$

\noindent\textbf{Lemma 1}. $\overline{\bar{R}_{1}\cup\bar{R}_{2}}%
=\overline{R_{1}\cup R_{2}}.$

\noindent Proof: We have that $R_{1}\cup R_{2}\subset\bar{R}_{1}\cup\bar
{R}_{2}$ so $\overline{R_{1}\cup R_{2}}\subset\overline{\bar{R}_{1}\cup\bar
{R}_{2}}$ while $\bar{R}_{1}\subset\overline{R_{1}\cup R_{2}},\bar{R}%
_{2}\subset\overline{R_{1}\cup R_{2}}$ implies $\overline{\bar{R}_{1}\cup
\bar{R}_{2}}\subset$ $\overline{R_{1}\cup R_{2}}.\smallskip$

\noindent This says that the equivalence relation generated by the union of
relations is equal to the equivalence relation generated by the union of the
corresponding generated equivalence relations. Furthermore, it is clear that
the union of equivalence relations is not in general an equivalence relation.

\section{Statistical Relations and Principles}

We define a \textit{statistical relation} to be a relation on $\mathcal{I}$
and a \textit{statistical principle} to be an equivalence relation on
$\mathcal{I}.$ The idea behind a statistical principle, as used here, is that
equivalent inference bases contain the same amount of statistical information
about the unknown $\theta.$ We make no attempt to give a precise definition of
what statistical information means. Birnbaum (1962) identified two inference
bases $I_{1},I_{2}\in\mathcal{I}$ as containing the same amount of statistical
information via the notation $Ev(I_{1})=Ev(I_{2}).$ We consider several
statistical relations.

The \textit{likelihood relation} $L$ on $\mathcal{I}$ is defined by
$(I_{1},I_{2})\in L$ whenever $\Theta_{E_{1}}=\Theta_{E_{2}}$ and there exists
$c>0$ such that $f_{E_{1},\theta}(x_{1})=cf_{E_{2},\theta}(x_{2})$ for every
$\theta.$ We have the following obvious result.$\smallskip$

\noindent\textbf{Lemma 2}. $L$ is a statistical principle.$\smallskip$

Actually the likelihood principle does not completely express the idea that
two inference bases with the same likelihood function contain the same amount
of statistical information. For this we need another statistical relation. We
define the \textit{invariance relation} $G$ by $(I_{1},I_{2})\in G$ whenever
there exist 1-1, onto, smooth functions $g:\mathcal{X}_{E_{1}}\rightarrow
\mathcal{X}_{E_{2}},h:\Theta_{E_{1}}\rightarrow\Theta_{E_{2}}$ with
$g(x_{1})=x_{2}$ and such that $f_{E_{1},\theta}(x)=f_{E_{2},h(\theta
)}(g(x))J_{g}^{-1}(x)$ for every $x\in\mathcal{X}_{E_{1}}$ where
$J_{g}(x)=(\det(\partial g(x)/\partial x))^{-1}=1$ in the discrete case. We
have the following result.$\smallskip$

\noindent\textbf{Lemma 3}. $G$ is a statistical principle.$\smallskip$

Now consider the equivalence relation $\overline{L\cup G}.$ If $(I_{1}%
,I_{2})\in L$ and $(I_{2},I_{3})\in G,$ then, for some constant $c>0$ and
mappings $g$ and $h,$ $f_{E_{1},\theta}(x_{1})=cf_{E_{2},\theta}%
(x_{2})=cf_{E_{3},h(\theta)}(g(x_{3}))J_{g}^{-1}(x_{3})=c^{\prime}%
f_{E_{3},h(\theta)}(g(x_{3}))$ and, so after relabelling, $I_{1}$ and $I_{3},$
have proportional likelihoods. Similarly, if $(I_{1},I_{2})\in G$ and
$(I_{2},I_{3})\in L,$ then again, after relabelling, $I_{1}$ and $I_{3}$ have
proportional likelihoods. So $(I_{1},I_{2})\in\overline{L\cup G}$ just
expresses the fact that $I_{1}$ and $I_{2}$ have proportional likelihoods,
perhaps after relabelling the data and the parameter. In this case we can
state clearly what the equivalence relation $\overline{L\cup G}$ expresses and
the generated equivalence relation makes sense. We do not need $\overline
{L\cup G},$ however, for a discussion of Birnbaum's result.

The \textit{sufficiency relation} $S$ is defined by $(I_{1},I_{2})\in S$
whenever $\Theta_{E_{1}}=\Theta_{E_{2}}$ and there exist minimal sufficient
statistics $m_{1}$ for $E_{1}$ and $m_{2}$ for $E_{2}$ such that the marginal
models induced by the $m_{i}$ are the same and $m_{1}(x_{1})=m_{2}(x_{2}).$ We
have the following result.$\smallskip$

\noindent\textbf{Lemma 4}. $S$ is a statistical principle and $S\subset L.$

\noindent Proof: Clearly\ $S$ is reflexive and symmetric and $S\subset L.$
Suppose $(I_{1},I_{2})\in S$ via the minimal sufficient statistics $m_{1}$ and
$m_{2}$ and $(I_{2},I_{3})\in S$ via the minimal sufficient statistics
$m_{2}^{\prime}$ and $m_{3}.$ Since any two minimal sufficient statistics are
1-1 functions of each other, there exists 1-1 function $h$ such that
$m_{2}^{\prime}=h\circ m_{2}.$ Then $(I_{1},I_{3})\in S$ via the minimal
sufficient statistics $h\circ m_{1}$ and $m_{3}.\smallskip$

\noindent Obviously we have the result that $(I_{1},I_{2})\in S$ whenever
$I_{2}$ can be obtained from $I_{1}$ via a sufficient statistic or conversely.
Furthermore, it makes sense to extend $S$ to $\overline{S\cup G}.$

The \textit{conditionality relation} $C$ is defined by $(I_{1},I_{2})\in C$
whenever $\Theta_{E_{1}}=\Theta_{E_{2}},$ $x_{1}=x_{2}$ and there exists
ancillary statistic $a$ for $E_{1}$ such that the conditional model given
$a(x_{1})$ is given by $E_{2}$ or with roles of $I_{1}$ and $I_{2}$ reversed.
We have the following result.$\smallskip$

\noindent\textbf{Lemma 5}. $C$ is reflexive and symmetric but is not
transitive and $C\subset L.$

\noindent Proof: The reflexivity, symmetry and $C\subset L$ are obvious. The
lack of transitivity follows via a simple example. Consider the model $E$ with
$\mathcal{X}_{E}=\{1,2\}^{2},\Theta_{E}=\{1,2\}$ and with $f_{E,\theta}$ given
by Table 1.
\begin{table}[h] \centering
\begin{tabular}
[c]{|c|c|c|c|c|}\hline
$(x_{1},x_{2})$ & $(1,1)$ & $(1,2)$ & $(2,1)$ & $(2,2)$\\\hline
$f_{E,1}(x_{1},x_{2})$ & $1/6$ & $1/6$ & $2/6$ & $2/6$\\\hline
$f_{E,2}(x_{1},x_{2})$ & $1/12$ & $3/12$ & $5/12$ & $3/12$\\\hline
\end{tabular}
\caption{Unconditional distributions.}%
\end{table}
Now note that $U(x_{1},x_{2})=x_{1}$ and $V(x_{1},x_{2})=x_{2}$ are both
ancillary and the conditional models, when we observe $(x_{1},x_{2})=(1,1),$
are given by Tables 2 and 3.
\begin{table}[h] \centering
\begin{tabular}
[c]{|c|c|c|c|c|}\hline
$(x_{1},x_{2})$ & $(1,1)$ & $(1,2)$ & $(2,1)$ & $(2,2)$\\\hline
$f_{E,1}(x_{1},x_{2}\,|\,U=1)$ & $1/2$ & $1/2$ & $0$ & $0$\\\hline
$f_{E,2}(x_{1},x_{2}\,|\,U=1)$ & $1/4$ & $3/4$ & $0$ & $0$\\\hline
\end{tabular}
\caption{Conditional distributions given $U=1$.}%
\end{table}
\begin{table}[h] \centering
\begin{tabular}
[c]{|c|c|c|c|c|}\hline
$(x_{1},x_{2})$ & $(1,1)$ & $(1,2)$ & $(2,1)$ & $(2,2)$\\\hline
$f_{E,1}(x_{1},x_{2}\,|\,V=1)$ & $1/3$ & $0$ & $2/3$ & $0$\\\hline
$f_{E,2}(x_{1},x_{2}\,|\,V=1)$ & $1/6$ & $0$ & $5/6$ & $0$\\\hline
\end{tabular}
\caption{Conditional distributions given $V=1$.}%
\end{table}%

\noindent The only ancillary for both these conditional models is the trivial
ancillary (the constant map). Therefore, there are no applications of $C$ that
lead to the inference base $I_{2},$ given by Table 2 with data $(1,1),$ being
related to the inference base $I_{3},$ given by Table 3 with data $(1,1).$ But
both of $I_{2}$ and $I_{3}$ are related under $C$ to the inference base
$I_{1}$ given by Table 1 with data $(1,1).$ This establishes the
result.$\smallskip$

\noindent Note that even under relabellings, the inferences bases $I_{2}$ and
$I_{3}$ in Lemma 5 are not equivalent.

If we are going to say that $(I_{1},I_{2})\in C$ means that $I_{1}$ and
$I_{2}$ contain an equivalent amount of information under $C,$ then we are
forced to expand $C$ to $\bar{C}$ so that it is an equivalence relation. But
this implies that the two inference bases $I_{2}$ and $I_{3}$ presented in the
proof of Lemma 5 contain an equivalent amount of information and yet they are
not directly related via $C.$ Rather they are related only because they are
conditional models obtained from a supermodel that has two essentially
different maximal ancillaries.

Saying that such models contain an equivalent amount of statistical
information is clearly a substantial generalization of $C.$ Note that, for the
example in the proof of Lemma 5, when $(1,1)$ is observed, the MLE is
$\hat{\theta}(1,1)=1.$ To measure the accuracy of this estimate we can compute
the conditional probabilities based on the two inference bases, namely,
\begin{align*}
P_{1}(\hat{\theta}(x_{1},x_{2}) &  =1\,|\,U=1)=1/2,P_{2}(\hat{\theta}%
(x_{1},x_{2})=2\,|\,U=1)=3/4\\
P_{1}(\hat{\theta}(x_{1},x_{2}) &  =1\,|\,V=1)=1/3,P_{2}(\hat{\theta}%
(x_{1},x_{2})=2\,|\,V=1)=5/6
\end{align*}
and so the accuracy of $\hat{\theta}$ is quite different depending on whether
we use $I_{2}$ or $I_{3}.$ It seems unlikely that we would interpret these
inference bases as containing an equivalent amount of information in a
frequentist formulation of statistics. As noted in Section 2, there is no
reason why we have to accept the equivalences given by a generated equivalence
relation unless we are certain that this equivalence relation expresses the
essence of the basic relation. It seems clear that there is a problem with the
assertion that $(I_{1},I_{2})\in\bar{C}$ means that $I_{1}$ and $I_{2}$
contain an equivalent amount of information without further justification.

We now follow a development similar to that found in Evans, Fraser and Monette
(1986) to prove the following result.$\smallskip$

\noindent\textbf{Theorem 6}. $C\subset\bar{C}=L$ where the first containment
is proper.

\noindent Proof: Clearly $C\subset\bar{C}$ and this containment is proper by
Lemma 5. If $(I_{1},I_{2})\in\bar{C},$ then (1) implies $(I_{1},I_{2})\in L$
since $C\subset L$ and so $\bar{C}\subset L.$ Now suppose that $(I_{1}%
,I_{2})\in L.$ We have that $f_{E_{1},\theta}(x_{1})=cf_{E_{2},\theta}(x_{2})$
for every $\theta$ for some $c>0.$ Assume first that $c>1.$ Now construct a
new inference base $I_{1}^{\ast}=(E_{1}^{\ast},(1,x_{1}))$ where
$\mathcal{X}_{E_{1}^{\ast}}=\{0,1\}\times\mathcal{X}_{E_{1}},$ and
$\{f_{E_{1}^{\ast},\theta}:\theta\in\Theta_{E_{1}}\}$ is given by Table 4
where $x_{10},x_{100},\ldots$ are the elements of $\mathcal{X}_{E_{1}}$ not
equal to $x_{1}$ and $p\in\lbrack0,1)$ satisfies $p/(1-p)=1/c.$
\begin{table}[h] \centering
\begin{tabular}
[c]{|c|c|c|c|c|}\hline
& $x_{1}$ & $x_{10}$ & $x_{100}$ & $\cdots$\\\hline
$i=1$ & $pf_{E_{1},\theta}(x_{1})$ & $pf_{E_{1},\theta}(x_{10})$ &
$pf_{E_{1},\theta}(x_{100})$ & $\cdots$\\\hline
$i=0$ & $1-p-pf_{E_{1},\theta}(x_{1})$ & $pf_{E_{1},\theta}(x_{1})$ & $0$ &
$\cdots$\\\hline
\end{tabular}
\caption{The model  $E_1^*$.}%
\end{table}
Then we see that $U(i,x)=i$ is ancillary as is $V$ given by $V(i,x)=1$ when
$x=x_{1}$ and $V(i,x)=0$ otherwise. Conditioning on $U(i,x)=1$ gives that
$(I_{1}^{\ast},I_{1})\in C$ while conditioning on $V(i,x)=1$ gives that
$(I_{1}^{\ast},I)\in C$ where $I=((\{0,1\},\{p_{\theta}:\theta\in\Theta
_{E_{1}}\}),1)$ and $p_{\theta}$ is the Bernoulli$(f_{E_{1},\theta}(x_{1})/c)$
probability function. Now, using $I_{2}$ we construct $I_{2}^{\ast}$ by
replacing $p$ by $1/2$ and $f_{E_{1},\theta}(x_{1})$ by $f_{E_{2},\theta
}(x_{2})$ in Table 4 and obtain that $(I_{2}^{\ast},I)\in C$ since
$f_{E_{1},\theta}(x_{1})/c=f_{E_{2},\theta}(x_{2}).$ Using (1) we have that
$(I_{1},I_{2})\in\bar{C}.$ If $c\leq1$ we start the construction process with
$I_{2}$ instead. This proves that $\bar{C}=L.\smallskip$

\noindent The proof that $L\subset\bar{C}$ relies on discreteness. This was
weakened in Evans, Fraser and Monette (1986) and even further weakened in Jang (2011).

We now show that Birnbaum's proof actually establishes the following
result.$\smallskip$

\noindent\textbf{Theorem 7}. $S\cup C\subset L\subset\overline{S\cup C}$

\noindent Proof: The first containment is obvious. For the second suppose that
$(I_{1},I_{2})\in L.$ We construct a new inference base $I=(E,y)$ from $I_{1}$
and $I_{2}$ as follows. Let $E$ be given by $\mathcal{X}_{E}=(1,\mathcal{X}%
_{E_{1}})\cup(2,\mathcal{X}_{E_{2}}),$%
\begin{align*}
f_{E,\theta}(1,x)  &  =\left\{
\begin{array}
[c]{cc}%
(1/2)f_{E_{1},\theta}(x) & \text{when }x\in\mathcal{X}_{E_{1}}\\
0 & \text{otherwise,}%
\end{array}
\right. \\
f_{E,\theta}(2,x)  &  =\left\{
\begin{array}
[c]{cc}%
(1/2)f_{E_{2},\theta}(x) & \text{when }x\in\mathcal{X}_{E_{2}}\\
0 & \text{otherwise.}%
\end{array}
\right.
\end{align*}
Then
\[
g(i,x)=\left\{
\begin{array}
[c]{cc}%
(i,x) & \text{when }x\notin\{x_{1},x_{2}\}\\
\{x_{1},x_{2}\} & \text{otherwise}%
\end{array}
\right.
\]
is sufficient for $E$ and so $((E,(1,x)),(E,(2,x)))\in S$ by the comment after
Lemma 4. Also, $h(i,x)=i$ is ancillary for $E$ and thus $((E,(1,x_{1}%
)),(E_{1},x_{1}))\in C$ and $((E,(2,x_{2})),(E_{2},x_{2}))\in C.$ Then by (1)
we have that $((E_{1},x_{1}),(E_{2},x_{2}))\in\overline{S\cup C}$ and we are
done.$\smallskip$

\noindent Note that Birnbaum's proof only proves the containments with no
equalities but we have the following result.$\smallskip$

\noindent\textbf{Theorem 8}. $S\cup C\ $is properly contained in $L$ while
$L=\overline{S\cup C}.$

\noindent Proof: To show that $S\cup C\subset L$ is proper, suppose that
$E_{1}$ is the Bernoulli$(\theta),\theta\in(0,1]$ model, $E_{2}$ is the
Geometric$(\theta),\theta\in(0,1]$ model and we observe $x_{1}=1$ and
$x_{2}=0$ so $f_{E_{1},\theta}(1)=\theta=f_{E_{2},\theta}(0).$ Note that the
full data is minimal sufficient for both $E_{1}$ and $E_{2}$ with
$\mathcal{X}_{E_{1}}=\{0,1\},\mathcal{X}_{E_{2}}=\{0,1,2,\ldots\}$ and further
that both of these models have only trivial ancillaries. Therefore, if
$I_{i}=(E_{i},1)$ we have that $(I_{1},I_{2})\notin S,(I_{1},I_{2})\notin C$
but $(I_{1},I_{2})\in L$ which proves that $S\cup C$ is properly contained in
$L.$

To prove that the second containment is exact we have, using (1), that
$(I_{1},I_{2})\in\overline{S\cup C}$ implies that $I_{1}$ and $I_{2}$ give
rise to proportional likelihoods as this is true for each element of $S\cup C$
and so $\overline{S\cup C}\subset L.\smallskip$

So we do not have, as usually stated for Birnabum's theorem, that $S$ and $C$
are together equivalent to $L$ but we do have that $\overline{S\cup C}$ is
equivalent to $L.$ Acceptance of $\overline{S\cup C}$ is not entailed,
however, by acceptance of both $S$ and $C$ as we have to examine the
additional relationships added to $S\cup C$ to see if they make sense. If one
wishes to say that acceptance of $S$ and $C$ implies the acceptance of
$\overline{S\cup C}$, then a compelling argument is required for these
additions and this seems unlikely. From the example of the proof of Theorem 8
we can see that acceptance of $\overline{S\cup C}$ is indeed equivalent to
acceptance of $L.$

From Theorems 6 and Theorem 7 we have the following Corollary.$\smallskip$

\noindent\textbf{Corollary 9}. $S\cup C\subset\bar{C}=L$ where the first
containment is proper. Furthermore, $S\subset\bar{C}$ and this containment is
proper.$\smallskip$

\noindent A direct proof that $S\subset\bar{C}$ has been derived by Jang
(2011). It is interesting to note that Corollary 9 shows that the existence of
$S$ in the modified statement of Birnbaum's theorem, where we require that we
accept all the equivalences generated by $S\ $\ and $C,$ is irrelevant as it
is not required. This is a reassuring result as it is unlikely that $S$ is
defective but it is almost certain that $C$ is defective, at least as
currently stated. Also we have the following result. $\smallskip$

\noindent\textbf{Lemma 10}. $\overline{C\cup G}=\overline{L\cup G}$

\noindent Proof: This is immediate from Lemma 1 and Lemma 3.$\smallskip$

\noindent This says that the equivalences obtained by combining invariance
under relabelling with conditionality are the same as the equivalences
obtained by combining invariance under relabelling with likelihood.

As with the proof of Birnbaum's theorem, the proof that $C=L$ provided in
Evans, Fraser and Monette (1986) is really a proof that $\bar{C}=L.$ This can
be seen from the proof of Theorem 6. So accepting the relation $C$ is not
really equivalent to accepting $L$ unless we agree that the additional
elements of $\bar{C}$ make sense. This is essentially equivalent to saying
that it doesn't matter which maximal ancillary we condition on and it is
unlikely that this is acceptable to most frequentist statisticians and this is
illustrated by the discussion concerning the example in Lemma 5.

As noted in Durbin (1970), requiring that any ancillaries used in an
application of $C$ be functions of a minimal sufficient statistic voids
Birnabum's proof, as the ancillary statistic used in the proof of Theorem 7 is
not a function of the sufficient statistic used in the proof. It is not clear,
however, what this restriction does to the result $\bar{C}=L,$ but we note
that there are situations where there exist nonunique maximal ancillaries
which are functions of the minimal sufficient statistic. In these
circumstances we would still be forced to conclude the equivalence of
inference bases derived by conditioning on the different maximal ancillaries
if we reasoned as in Evans, Fraser and Monette (1986). Of course, we are
arguing here that the result requires the statement of an additional hypothesis.

\section{Conclusions}

We have shown that the proof in Birnbaum (1962) did not prove that $S$ and $C$
lead to $L.$ Rather the proof establishes that $\overline{S\cup C}=L$ and this
is something quite different. The statement of Birnbaum's theorem in prose
should have been: if we accept the relation $S$ and we accept the relation $C$
\textit{and} we accept all the equivalences generated by $S$ and $C$ together,
then this is equivalent to accepting $L.$ The essential flaw in Birnbaum's
theorem lies in excluding this last hypothesis from the statement of the
theorem. The same qualification applies to the result proved in Evans, Fraser
and Monette (1986) where the statement of the theorem should have been: if we
accept the relation $C$ \textit{and} we accept all the equivalences generated
by $C,$ then this is equivalent to accepting $L.$

The way out of the difficulties posed by Birnbaum's theorem, and the result
relating $C$ and $L,$ is to acknowledge that additional hypotheses are
required for the results to hold. Certainly these results seem to lose their
impact when they are correctly stated and we realize that an equivalence
relation generated by a relation is not necessarily meaningful. It is
necessary to provide an argument as to why the generated equivalence relation
captures the essence of the relation that generates it and it is not at all
clear how to do this in these cases.

As we have noted, the essential result in all of this is $\bar{C}=L$ and this
has some content albeit somewhat minor. Furthermore, the proof of this result
is based on a defect in $C,$ namely, it is not an equivalence relation due to
the general nonexistence of unique maximal ancillaries. As such it is hard to
accept $C$ as stated as any kind of characterization of statistical
evidence.\ Given the intuitive appeal of this relation in some simple
examples, however, resolving the difficulties with $C$ still poses a major
challenge for a frequentitst theory of statistics.

\begin{center}
\textbf{{\Large References}}\bigskip
\end{center}

\noindent Barndorff-Nielsen, O. E. (1995) Diversity of evidence and Birnbaum's
theorem (with discussion). Scand. J. Statist. 22 (4), 513--522.\medskip

\noindent Birnbaum, A. (1962) On the foundations of statistical inference
(with discussion). J. Amer. Stat. Assoc., 57, 269-332.\medskip

\noindent Durbin, J. (1970) On Birnbaum's theorem on the relation between
sufficiency, conditionality and likelihood. J. Amer. Stat. Assoc., 654,
395-398.\medskip

\noindent Evans, M., Fraser, D.A.S. and Monette, G. (1986) On principles and
arguments to likelihood (with discussion). Canad. J. of Statsitics, 14, 3,
181-199.\medskip

\noindent Gandenberger, G. (2012) A new proof of the likelihood principle. To
appear in British Journal for Philosophy of Science.\medskip

\noindent Helland, I.S. (1995) Simple counterexamples against the
conditionality principle. Amer. Statist., 49, 4, 351-356.\medskip

\noindent Holm, S. (1985) Implication and equivalence among statistical
inference rules. In Contributions to Probability and Statistics in Honour of
Gunnar Blom. Univ. Lund, Lund, 143--155.\medskip

\noindent Jang, G. H. (2011) The conditionality principle implies the
sufficiency principle. Working paper.\medskip

\noindent Kalbfleisch, J.D. \ (1975) Sufficiency and conditionality.
Biometrika, 62, 251-259.\medskip

\noindent Mayo, D. (2010). An Error in the Argument from Conditionality and
Sufficiency to the Likelihood Principle. In Error and Inference: Recent
Exchanges on Experimental Reasoning, Reliability and the Objectivity and
Rationality of Science (D. Mayo and A. Spanos eds.), Cambridge: Cambridge
University Press: 305-14.

\end{document}